\documentclass[a4paper,12pt]{article}
    \usepackage[top=2.5cm,bottom=2.5cm,left=2.5cm,right=2.5cm]{geometry}
\usepackage{graphicx}
\usepackage{latexsym,bm}
\usepackage{amsthm,amsmath,amssymb,datetime}
\begin{document}

\newtheorem{prop}{Proposition}[section]
\renewcommand{\theprop}{\arabic{section}.\arabic{prop}}
\newtheorem{lem}[prop]{Lemma}
\renewcommand{\thelem}{\arabic{section}.\arabic{lem}}
\newtheorem{thm}[prop]{Theorem}
\renewcommand{\thethm}{\arabic{section}.\arabic{thm}}
\newtheorem{cor}[prop]{Corollary}
\renewcommand{\thecor}{\arabic{section}.\arabic{cor}}
\renewcommand{\theequation}{\arabic{section}.\arabic{equation}}
\renewcommand{\thefigure}{\arabic{figure}}
\newtheorem{exa}{prop}
\renewcommand{\theexa}{\arabic{section}.\arabic{ex}}
\newtheorem{rem}[prop]{Remark}
\renewcommand{\therem}{\arabic{section}.\arabic{rem}}
\renewcommand{\theequation}{
\arabic{equation}}

\newtheorem{defi}[prop]{Definition}
\renewcommand{\thedefi}{\arabic{section}.\arabic{defi}}
\renewcommand{\thefigure}{\arabic{figure}}
\def\zm{\noindent{\bf Proof.\ }}
\def\ezm{\vspace*{6mm}\framebox{} }
\renewcommand{\deg}{\overline{d}} 

\addtolength{\baselineskip}{+0.6mm}

\title{{\bf On the connective eccentricity index of two types of trees}\thanks{Corresponding author: hydeng@hunnu.edu.cn (Hanyuan Deng).
Project supported by the program for excellent talents in
Hunan Normal University(ET13101) and the National Natural Science Foundation of China (61572190).}}
\author{Zikai Tang, Lingyao Jiang, Hanyuan Deng*
\\{\small College of Mathematics and Computer Science,}
\\{\small Hunan Normal University, Changsha, Hunan 410081, P. R.
China}}

\date{}
\maketitle

{\bf Abstract}: The connective eccentricity index  $\xi^{ce}=\sum^{}_{u\in V}\frac{d(u)}{\varepsilon(u)}$, where
 $\varepsilon(u)$ and $d(u)$ denote the eccentricity and the degree of the vertex $u$, respectively. In this paper, we first determine the extremal trees which minimize and maximize the connective eccentricity index among all trees with a given degree sequence, and then determine the extremal trees which minimize and maximize the connective eccentricity index among all trees with a given number of branching vertices.


\addtolength{\baselineskip}{+1mm}

\section{Introduction}

Let $G$ be a simple and connected graph with $n=|V(G)|$ vertices.
For a vertex $v\in V(G)$, $d_G(v)$ denotes the degree of $v$ (or just $d(v)$ briefly). For
vertices $v, u\in V(G)$, the distance $d(v,u)$ is defined as the
length of a shortest path between $v$ and $u$ in $G$. The
eccentricity $\varepsilon_G(v)$ (or just $\varepsilon(v)$ briefly) of a vertex $v$ is the maximum
distance from $v$ to any other vertex of $G$. The diameter $D(G)$ of a graph
is the maximum eccentricity of any vertex in the graph. A vertex of degree one is
called a pendant vertex. A path $P=v_0v_1\cdots v_t$ of $G$ is a
pendant path if $v_0$ is a pendant vertex, the degree of any internal
vertex $v_i$ ($1\leq i<t$) is two and the degree of $v_t$ is at
least three. Let $S_n$ and $P_n$ denote the star and the path with
$n$ vertices, respectively. For other terminologies and notations
not defined here we refer the readers to \cite{bm76}.

In 2000, Gupta, Singh and Madan \cite{SMA2000} introduced a novel, adjacency-cum-path length
based, topological descriptor termed the connective eccentricity index. In order to explore
the potential of the connective eccentricity index in predicting biological activity, authors
used nonpeptide N-benzylimidazole derivatives to investigate the predictability of the
connective eccentricity index with respect to antihypertensive activity. They showed that
results obtained using the connective eccentricity index were better than the corresponding
values obtained using Balaban's mean square distance index and the accuracy of prediction
was found to be about $80\%$ in the active range \cite{SMA2000}.

The connective eccentricity index (CEI) of a graph G was defined as
$$\xi^{ce}(G)=\sum^{}_{u\in V}\frac{d(u)}{\varepsilon(u)}=\sum^{}_{uv\in E}(\frac{1}{\varepsilon(u)}+\frac{1}{\varepsilon(v)}).\eqno(1)$$

The upper or lower bounds for the connective eccentricity index in terms of some graph invariants such as the radius, the independence number, the vertex connectivity, the minimum degree, the maximum degree etc. were recently reported in \cite{FY13,YQT14,XDL16}. In this paper, we will prove that the "greedy" caterpillar minimizes $\xi^{ce}(T)$, while the "greedy" tree maximizes $\xi^{ce}(T)$ among all trees with a given degree sequence. Moreover, we will determine the lower and upper bounds for the connective eccentricity index of an $n$-vertex tree with a given number of branching vertices.

\section{ Preliminaries}

In the following, we give some transformations which will be used in the next section.

\begin{lem}\label{lem2.1} ({\bf The transformation A})
Let $u$ be a vertex of a graph $Q$ with at least two vertices. For
integer $a\geq 1$, $G_1$ is the tree obtained by attaching a star
$S_{a+1}$ at its center $v$ to $u$ of $Q$, $G_2$ is the tree obtained
by attaching $a+1$ pendent vertices to $u$ of $Q$ (see Figure 1). Then $\xi^{ce}(G_2)\geq \xi^{ce}(G_1)$.
\end{lem}
\zm  In the graph $G_1$, we have $\varepsilon_{G_1}(v)=\varepsilon_{Q}(u)+1\geq \varepsilon_{G_1}(u)=max\{2,\varepsilon_{Q}(u)\}$. It is easy to see from Figure 1 that $\varepsilon_{G_1}(w)\geq \varepsilon_{G_2}(w)$ for $w\in V$ and $d_{G_1}(w)=d_{G_2}(w)$ for $w\in V-\{u,v\})$.
By the definition of $\xi^{ce}(G)$, we have
\[\begin{array}{ll}
\xi^{ce}(G_2)-\xi^{ce}(G_1)&\geq \frac{d_{G_2}(u)}{\varepsilon_{G_2}(u)}-\frac{d_{G_1}(u)}{\varepsilon_{G_1}(u)}+\frac{d_{G_2}(v)}{\varepsilon_{G_2}(v)}
-\frac{d_{G_1}(v)}{\varepsilon_{G_1}(v)}\\
&\geq\frac{d_{G_2}(u)}{\varepsilon_{G_1}(u)}-\frac{d_{G_1}(u)}{\varepsilon_{G_1}(u)}+\frac{d_{G_2}(v)}{\varepsilon_{G_1}(v)}
-\frac{d_{G_1}(v)}{\varepsilon_{G_1}(v)}\\
&=a(\frac{1}{\varepsilon_{G_1}(u)}-\frac{1}{\varepsilon_{G_1}(v)})\\
&\geq 0.
\end{array}\]
\ezm

\begin{lem}\label{lem2.2} ({\bf The transformation B})
Let $w$ be a vertex of a nontrivial connected graph $G$. For
nonnegative integers $p$ and $q$, $G(p, q)$ denotes the graph
obtained from $G$ by attaching to the vertex $w$ pendent paths $P=
wv_1v_2 \cdots v_p$ and $Q=wu_1u_2 \cdots u_q$ of lengths $p$ and
$q$, respectively. If $p\geq q\geq 1$, then $\xi^{ce}(G(p,q))
\geq\xi^{ce}(G(p+1,q-1))$.
\end{lem}
\zm (1) For $q>1$, let $G(p+1,q-1)$ be obtained from $G(p,q)$ by deleting the edge $u_{q}u_{q-1}$ and adding an edge $v_{p}u_{q}$. We have $\varepsilon_{G(p,q)}(v_{p})\geq \varepsilon_{G(p+1,q-1)}(v_{p})$, $\varepsilon_{G(p,q)}(u_{q-1})\leq \varepsilon_{G(p+1,q-1)}(u_{q-1})$ and $\varepsilon_{G(p,q)}(u_{q-1})\leq \varepsilon_{G(p+1,q-1)}(v_{p})$.  If $t\in V-\{u_p,v_{q-1}\}$, then $d_{G(p,q)}(t)=d_{G(p+1,q-1)}(t)$ and $\varepsilon_{G(p,q)}(t)\leq \varepsilon_{G(p+1,q-1)}(t)$. So,
\[
\begin{array}{l}
\xi^{ce}(G(p+1,q-1))-\xi^{ce}(G(p,q)\\
\leq \frac{d_{G(p+1,q-1)}(v_p)}{\varepsilon_{G(p+1,q-1)}(v_p)}-\frac{d_{G(p,q)}(v_p)}{\varepsilon_{G(p,q)}(v_p)}+
\frac{d_{G(p+1,q-1)}(u_{q-1})}{\varepsilon_{G(p+1,q-1)}(u_{q-1})}
-\frac{d_{G(p,q)}(u_{q-1})}{\varepsilon_{G(p,q)}(u_{q-1})}+
\frac{d_{G(p+1,q-1)}(u_{q})}{\varepsilon_{G(p+1,q-1)}(u_{q})}
-\frac{d_{G(p,q)}(u_{q})}{\varepsilon_{G(p,q)}(u_{q})}\\
=\frac{2}{\varepsilon_{G(p+1,q-1)}(v_p)}-\frac{1}{\varepsilon_{G(p,q)}(v_p)}+
\frac{1}{\varepsilon_{G(p+1,q-1)}(u_{q-1})}-\frac{2}{\varepsilon_{G(p,q)}(u_{q-1})}+
\frac{1}{\varepsilon_{G(p+1,q-1)}(u_{q})}
-\frac{1}{\varepsilon_{G(p,q)}(u_{q})}\\
\leq\frac{3}{\varepsilon_{G(p+1,q-1)}(u_{q})}-\frac{3}{\varepsilon_{G(p,q)}(v_p)}\\
\leq 0.
\end{array}\]

(2) For $q=1$, let $G(p+1,0)$ be obtained from $G(p,1)$  by deleting the edge $u_{1}w$ and adding an edge $v_{p}u_{1}$. If $t\in V-\{u_p,w\}$, then $d_{G(p,q)}(t)=d_{G(p+1,q-1)}(t)$ and $\varepsilon_{G(p,q)}(t)\leq \varepsilon_{G(p+1,q-1)}(t)$. So, we have
\[
\begin{array}{l}
\xi^{ce}(G(p+1,0))-\xi^{ce}(G(p,1)\\
\leq \frac{d_{G(p+1,0)}(v_p)}{\varepsilon_{G(p+1,0)}(v_p)}-\frac{d_{G(p,1)}(v_p)}{\varepsilon_{G(p,1)}(v_p)}+
\frac{d_{G(p+1,0)}(w)}{\varepsilon_{G(p+1,0)}(w)}
-\frac{d_{G(p,1)}(w)}{\varepsilon_{G(p,1)}(w)}+
\frac{d_{G(p+1,0)}(u_{1})}{\varepsilon_{G(p+1,0)}(u_{1})}
-\frac{d_{G(p,1)}(u_{1})}{\varepsilon_{G(p,1)}(u_{1})}\\
=\frac{2}{\varepsilon_{G(p+1,0)}(v_p)}-\frac{1}{\varepsilon_{G(p,1)}(v_p)}+
\frac{d(w)-1}{\varepsilon_{G(p+1,0)}(w)}-\frac{d(w)}{\varepsilon_{G(p,1)}(w)}+
\frac{1}{\varepsilon_{G(p+1,0)}(u_{1})}
-\frac{1}{\varepsilon_{G(p,1)}(u_{1})}\\
\leq\frac{3}{\varepsilon_{G(p+1,0)}(u_{1})}-\frac{3}{\varepsilon_{G(p,q)}(v_p)}\\
\leq 0.
\end{array}\]

From above, the result is proved.\ezm

By using Lemma \ref{lem2.1} and Lemma \ref{lem2.2}, we can obtain the following result directly.
\begin{prop}\label{T-second}
Let $T$ be a tree with $n\geq 6$ vertices and $T\neq S_n, P_n, T_1, T_2$ (depicted in Figure 2).  Then
$$\xi^{ce}(S_n)>\xi^{ce}(T_2)> \xi^{ce}(T)>\xi^{ce}(T_1)>\xi^{ce}(P_n).$$
\end{prop}

\section{The connective eccentricity index of trees with a given degree sequence}

Given a degree sequence, let $\mathcal{T}$ be the class of trees that realize this degree sequence. We will determine the trees which maximize or minimize the connective eccentricity index in $\mathcal{T}$, and will compare the maximal values of the connective eccentricity index for different degree sequences. Note that a sequence $(d_1, d_2,\cdots, d_n)$ of positive integers is a degree sequence of a tree if and only if $\sum_{i=1}^{n}d_i=2(n-1)$.

In the following, we firstly show that the greedy caterpillar minimize the connective eccentricity index in $\mathcal{T}$.

In \cite{WH2014}, Wang gave the definition of the {\it greedy caterpillar}. Greedy caterpillars are not unique with given a degree sequence.

\begin{defi}\cite{WH2014}
For $n\geq 3$, let $\overline{d} =(d_1, d_2, \ldots, d_n)$ be the non-increasing degree sequence of a tree with $d_k>1$ and $d_{k+1}=1$ for some $k\in\{1,2,\cdots,n-2\}$. The {\it greedy caterpillar}, $T$, is constructed as follows:
\begin{itemize}
 \item Start with a path $P=z_1z_2\ldots z_k$.
 \item Let $\phi: \{z_i\}_{i=1}^k \rightarrow \{d_i\}_{i=1}^k$ be a one-to-one function such that, for each pair $i,j\in [k]$, if $\varepsilon_P(z_i)> \varepsilon_P(z_{j})$, then $\phi(z_i)\geq \phi(z_{j})$ .
  \item For each $i\in \{2,3,\ldots, k-1\}$, attach $\phi(z_i)-2$ pendant vertices to $z_i$. For $i\in \{1,k\}$, attach $\phi(z_i)-1$ pendant vertices to $z_i$.
\end{itemize}
 \label{gr_cat}
\end{defi}

\begin{thm}\label{thm3-2}
Among trees with a given tree degree sequence, the greedy caterpillar has the minimum the connective eccentricity index.
\end{thm}
\zm Fix a degree sequence $\overline{d}=(d_1, \ldots, d_n)$ which is written in the form described in Definition \ref{gr_cat}. Let $\mathcal{T}$ be the collection of trees with degree sequence $\overline{d}$, and $T\in\mathcal{T}$ such that $\xi^{ce}(T)=\min_{F\in\mathcal{T}} \xi^{ce}(F)$. We first show that $T$ is a caterpillar.

By contradiction, suppose $T$ is not a caterpillar. Let $P_T(u,v)=uu_1u_2 \ldots u_{k}v$ be a longest path in $T$. Let $x\in\{1,2,\cdots,k\}$ be the least integer such that $u_x$ has a non-leaf neighbor $w$ not on $P_T(u,v)$. Then $x\neq 1$ for the maximality of $P_T(u,v)$. Let $W$ be the component containing $w$ in $T-\{u_xw\}$.

Create a new tree $T'$ from $T$ by replacing each edge of the form $zw$ in $W$ with the edge $zu$ (see Figure 3). Notice that $T$ and $T'$ have the same degree sequence. However, for any vertex $s \in \left(V(T) \setminus V(W)\right) \cup \{w\}$, $\varepsilon_{T'}(s) \geq \varepsilon_{T}(s)$ since $P_T(u,v)$ is a longest path in $T$. For any vertex $r\in V(W)-w$, we have
$$\varepsilon_{T'}(r) = d(r,u) + d(u,v) > d(u,v) \geq \varepsilon_T(r) . $$
By the definition of the connective eccentricity index, we have $ \xi^{ce}(T') <\xi^{ce}(T)$, a contradiction.

Now, we will show that $T$ is a greedy caterpillar. By contradiction, suppose $T$ is not a greedy caterpillar. Since $T$ is a caterpillar with internal vertices forming path $P=u_1u_2\ldots u_k$, the eccentricity of any internal vertex is independent of the interval vertex degree assignments. There must be $i,j\in\{1,2,\cdots,k\}$ with $d_T(u_i)<d_T(u_j)$ and $\varepsilon_T(u_i)>\varepsilon_T(u_j)$.

Create a new tree $T''$ from $T$ by  replacing each edge of the form $u_jw, u_it$ ($w, t$ be the pendant vertices of $u_j, u_i$, respectively) with the edge $u_iw, u_jt$. Notice that $T$ and $T''$ have the same degree sequence and $d_{T''}(u_i)=d_T(u_j)$, $d_{T''}(u_j)=d_T(u_i)$, $\varepsilon_{T''}(u_i)=\varepsilon_T(u_i)$, $\varepsilon_{T''}(u_j)=\varepsilon_T(u_j)$. We have

\[
\begin{array}{ll}\xi^{ce}(T'') -\xi^{ce}(T)&< \frac{d_{T''}(u_i)}{\varepsilon_{T''}(u_i)}+\frac{d_{T''}(u_j)}{\varepsilon_{T''}(u_j)}-
\frac{d_{T}(u_i)}{\varepsilon_{T}(u_i)}-\frac{d_{T}(u_j)}{\varepsilon_{T}(u_j)}\\
&=(d_{T}(u_j)-d_{T}(u_i))(\frac{1}{\varepsilon_{T}(u_i)}-\frac{1}{\varepsilon_{T}(u_j)})\\
&<0
\end{array}
\]
a contradiction. \ezm

Next, we will show that the greedy tree maximize the connective eccentricity index in $\mathcal{T}$.

Each tree is rooted at a vertex (while the root has no bearing on the connective eccentricity index, we use the added structure to direct our conversation).
The height of a vertex is the distance to the root, and the tree's height, $h=h(T)$, is the maximum of all heights of vertices. We start with some definitions.

\begin{defi}\cite{ssw14}\label{def:leveldegree} In a rooted tree, the list of multisets $L_i$ of degrees of vertices
at height $i$, starting with $L_0$ containing the degree of the root vertex, is called the {\it level-degree sequence} of the rooted tree.
\end{defi}

Let $|L_i|$ be the number of entries in $L_i$. It is easy to see that a list of multisets is the level degree sequence of a rooted tree if and only if (1) the multiset $\bigcup_i L_i$ is a tree degree sequence; (2) $|L_0|=1$; and (3)
$\sum_{d\in L_0} d=|L_{1}|$ and for all $i\geq 1$, $\sum_{d\in L_i} (d-1)=|L_{i+1}|$.

In a rooted tree, the \textit{down-degree} of the root is equal to its degree. The down degree of any other vertex is its degree minus one.

\begin{defi}\cite{ssw14}\label{lecel-greedy-tree}
Given the level-degree sequence of a rooted tree, the {\it level-greedy rooted tree} for this level-degree sequence is built as follows:
(1) For each $i\in\{1,2,\cdots,n\}$, place $|L_i|$ vertices in level $i$ and to each vertex, from left to right, assign a degree from $L_i$ in non-increasing order; (2) For $i\in\{1,2,\cdots,n-1\}$, from left to right, join the next vertex in $L_i$ whose down-degree is $d$ to the first $d$ so far unconnected vertices on level $L_{i+1}$. Repeat for $i+1$.
\end{defi}

\begin{defi}\cite{ssw14}\label{def_greedy}
Given a tree degree sequence $(d_1, d_2,\cdots, d_n)$ in non-increasing order, the {\it greedy tree} for this degree sequence is the level-greedy tree for the level-degree sequence that has $L_0=\{d_1\}$, $L_1=\{d_2,\cdots, d_{d_1+1}\}$ and for each $i>1$,
$$|L_i|=\sum_{d\in L_{i-1}} (d-1)$$
with every entry in $L_{i}$ at most as large as every entry in $L_{i-1}$.
\end{defi}

A greedy tree with the degree sequence $(4, 4, 4, 3, 3, 3, 3, 3, 3, 3, 2, 2, 1, \cdots, 1)$ is shown in Figure 4.

\begin{lem}\label{level greedy tree}
Among all the trees with a given level-degree sequence, the level-greedy tree maximizes the connective eccentricity index.
\end{lem}

\zm By induction on the number of vertices, the base case with one vertex is trivial.

Let $T$ be a rooted tree with the given level-degree sequence and maximize the connective eccentricity index (i.e. $T$ is optimal).
For vertices $w\in T_1$ and $u\in T-T_1$, both of height $j$ (See Figure 5), we notice that $\varepsilon_T(u)=j+h$, $\varepsilon_T(w)= \max\{j+h', \varepsilon_{T_1}(w)\}\leq \varepsilon_T(u)$.  Suppose for contradiction that $d_T(u)>d_T(w)$. Create a new tree $T'$ by moving $d_T(u)-d_T(w)$ children of $u$ and their descendants to adoptive parent $w$. This effectively switches the degrees of $u$ and $w$ while maintaining the level degree sequence.

While $\varepsilon_{T'}(u) = \varepsilon_T(u)$, notice that $h'$ does not increase and $\varepsilon_T(x)\geq \varepsilon_{T'}(x)$ for all $x \in V$. Since $\varepsilon_{T'}(w) \leq \max\{j+h', \varepsilon_{T_1}(w)\} =\varepsilon_T(w)$,  we have
 \[\begin{array}{ll}
 \xi^{ce}(T')-\xi^{ce}(T)&\geq \frac{d_{T'}(u)}{\varepsilon_{T'}(u)}+\frac{d_{T'}(w)}{\varepsilon_{T'}(w)}
 -\frac{d_{T}(u)}{\varepsilon_{T}(u)}-\frac{d_{T}(w)}{\varepsilon_{T}(w)}\\
 &\geq (d_T(w)-d_T(u))(\frac{1}{\varepsilon_T(u)}-\frac{1}{\varepsilon_T(w)})\\
 &>0
 \end{array}
 \]
a contradiction to the optimality of $T$. Otherwise, $T'$ and $T$ are both optimal trees. In this case, we can repeat this shifting of degrees for pairs of vertices of height 1, followed by pairs of vertices of height 2, and so on until we either meet a contradiction or construct an optimal tree in which $d(u) \leq d(w)$ for all $w\in T_1$ and $u\in T-T_1$ of the same height.

Now, we have a partition of the level-degree sequence for $T$ into the level-degree sequences for $T-T_1$. By the inductive hypothesis, we may assume that both  $T_1$ and $T-T_1$ are level-greedy trees on their level-degree sequences. As a result, $T$ is a level-greedy tree.
\ezm

The next theorem also yields a stronger result than merely the connective eccentricity  index among trees with a given degree sequence.

\begin{thm}\label{greedy-t}
Among all trees with a given degree sequence, the greedy tree has the maximal connective eccentricity index.
\end{thm}
\zm Let $\bar{d}=(d_1, d_2, \ldots, d_n)$ be given degree sequence in non-increasing order and $T^{*}$ the tree with the maximal connective eccentricity index with the given degree sequence.

Take a longest path in $T^*$ and a center vertex $v$ of this path as the root of $T^*$. In $T^*-\{v\}$, let $T_1$ be the component containing the leaf with the greatest height. By our choice of the root, if $h$ is the height of $T_1$, then $T-T_1$ has height $h'\in \{h-1,h\}$.  The vertex set $V$ of $T^{*}$ can be divided into $h$ subsets $V=V_0\cup V_1\cdots \cup V_{h}$, where $V_0=\{v\}$, $V_1=\{u_1,\cdots,u_{d(v)}\}$ and for each $i>1$,
\[|V_i|=\sum_{u\in V_{i-1}} d(u)-|V_{i-1}|.\]
By Lemma \ref{level greedy tree}, $T^{*}$ is a level greedy tree. Next, we will prove that degree of every entry in $V_{i}$ at most as large as degree of every entry in $V_{i-1}$.

Suppose that there are $V_i=\{w_1,\cdots,w_k\}$ and $V_{i-1}=\{v_1,\cdots,v_t\}$ such that $d_{T^*}(w_1)>d_{T^*}(v_t)$ and $w_1\in T_1$.  Create a new tree $T'$ by moving $d_{T^*}(w_1)-d_{T^*}(v_t)$ children of $w_1$ and their descendants to adoptive parent $v_t$ with the height of $T_1$ no change. This effectively switches the degrees of $w_1$ and $v_t$ while maintaining the degree sequence. We now examine two cases: $v_t\in T_1$ and $v_t\in T^*- T_1$.

{\bf Case I}. $v_t\in T_1$. Note that $\varepsilon_{T^*}(w_1)=\varepsilon_{T'}(w_1)\geq\varepsilon_{T^*}(v_t)=\varepsilon_{T'}(v_t)$ and $\varepsilon_{T'}(x)\leq \varepsilon_{T^*}(x)$ for all $x \in V$, we have
 \[\begin{array}{ll}
 \xi^{ce}(T')-\xi^{ce}(T^*)&\geq \frac{d_{T'}(v_t)}{\varepsilon_{T'}(v_t)}+\frac{d_{T'}(w_1)}{\varepsilon_{T'}(w_1)}
 -\frac{d_{T^*}(v_t)}{\varepsilon_{T^*}(v_t)}-\frac{d_{T^*}(w_1)}{\varepsilon_{T^*}(w_1)}\\
 &\geq (d_{T^*}(w_1)-d_{T^*}(v_1))(\frac{1}{\varepsilon_{T^*}(v_t)}-\frac{1}{\varepsilon_{T^*}(w_1)})\\
 &>0
 \end{array}
 \]
a contradiction to the optimality of $T$.

{\bf Case II}. $v_t\in T^*- T_1$. If $h'=h$, we notice that $\varepsilon_{T'}(w_1)=\varepsilon_{T^*}(w_1)=i+h> \varepsilon_{T'}(v_t)=\varepsilon_{T^*}(v_t)=i-1+h$ and $\varepsilon_{T'}(x)\leq \varepsilon_{T^*}(x)$ for all $x \in V$, then
\[\begin{array}{ll}
 \xi^{ce}(T')-\xi^{ce}(T^*)&\geq \frac{d_{T'}(v_t)}{\varepsilon_{T'}(v_t)}+\frac{d_{T'}(w_1)}{\varepsilon_{T'}(w_1)}
 -\frac{d_{T^*}(v_t)}{\varepsilon_{T^*}(v_t)}-\frac{d_{T^*}(w_1)}{\varepsilon_{T^*}(w_1)}\\
 &\geq (d_{T^*}(w_1)-d_{T^*}(v_1))(\frac{1}{\varepsilon_{T^*}(v_t)}-\frac{1}{\varepsilon_{T^*}(w_1)})\\
 &>0
 \end{array}
 \]
a contradiction to optimality of $T$.

If $h'=h-1$, we notice that $\varepsilon_{T'}(w_1)=\varepsilon_{T^*}(w_1)= \varepsilon_{T'}(v_t)=\varepsilon_{T^*}(v_t)=i-1+h$ and $\varepsilon_{T'}(x)= \varepsilon_{T^*}(x)$ for all $x \in V$, then $\xi^{ce}(T')=\xi^{ce}(T^*)$.

In conclusion, we have that the greedy tree has the maximal connective eccentricity index among the trees with a given degree sequence.\ezm

\begin{rem}
Such extremal trees are not necessarily unique. In fact, the greedy tree give a more stronger restriction than what we needed, as stated in the theorem, while still not being the unique structure.
\end{rem}

In the following, we will compare the connective eccentricity indices of greedy trees with different degree sequences.

\begin{defi}
 Let $\pi'=(d'_1, \cdots d'_{n})$ and $\pi''=(d''_1, \cdots,
 d''_{n})$ be two non-increasing tree degree sequences. $\pi''$ is said to {\it majorize} $\pi'$, denoted $\pi' \vartriangleleft \pi'' $, if for $k\in\{1,2,\cdots,n-1\}$
\[
  \sum_{i=0}^{k}d'_i\le\sum_{i=0}^k d''_i  \text{and~} \sum_{i=0}^{n}d'_i=\sum_{i=0}^{n}d''_i. \]
\end{defi}

\begin{lem}\cite{wei82}\label{lem3-3}
Let $\pi'=(d'_1, \cdots d'_{n})$ and $\pi''=(d''_1, \cdots,
d''_{n})$ be two non-increasing tree degree sequences. If
$\pi'\triangleleft \pi''$, then there exists a series of
(non-increasing) tree degree sequences $\pi^{(i)}=(d_1^{(i)}, \cdots, d_{n}^{(i)})$ for $1\leq i\leq m$ such that
 \[\pi' =\pi^{(1)}\triangleleft \pi^{(2)} \triangleleft \cdots \triangleleft \pi^{(m-1)} \triangleleft\pi^{(m)}= \pi''.\]
In addition, each $\pi^{(i)}$ and $\pi^{(i+1)}$ differ at exactly two entries, say the $j$ and $k$ entries, $j<k$, where $d_j^{(i+1)} = d_j^{(i)} +1$ and $d_k^{(i+1)} = d_k^{(i)} -1$.
\end{lem}

\begin{thm}\label{thm3-3}
Let $\pi'=(d'_1, \cdots d'_{n})$ and $\pi''=(d''_1, \cdots,
d''_{n})$ be two non-increasing greedy tree degree sequences. If
$\pi'\triangleleft \pi''$, then
\[\xi^{ce}(T_{\pi'}^*)\leq \xi^{ce}(T_{\pi''}^*)\]
where $T_{\pi}^*$ is the greedy tree for degree sequence $\pi$.
\end{thm}
\zm According to Lemma~\ref{lem3-3}, it suffices to compare the connective eccentricity indices of two greedy trees whose degree sequences differ in two entries, each by exactly 1, i.e., we can assume that
$$ \pi'=(d'_1, \cdots d'_{n}) \triangleleft (d_1'', \cdots, d_{n}'')=\pi'' $$
with $d''_j = d'_j +1$, $d''_k = d'_k -1$ for some $j<k$ and all other entries are the same.

Let $u$ and $v$ be the vertices corresponding to $d'_j$ and $d'_k$, respectively, and $w$ be a child of $v$ in $T_{\pi'}^*$ (see Figure 6).
Construct $T_{\pi''}$ from $T_{\pi'}^*$ by removing the edge $vw$ and adding edge $uw$. Note that $T_{\pi''}$ has the degree sequence $\pi''$, and by Theorem~\ref{greedy-t}
$$ \xi^{ce}(T_{\pi''}^*) \geq \xi^{ce}(T_{\pi''}) . $$

On the other hand, from the definition of the connective eccentricity index, we have
$$\begin{array}{ll}
\xi^{ce}(T_{\pi''}) - \xi^{ce}(T_{\pi'}^*) &\geq \frac{d_v''}{\varepsilon_{T_{\pi''}}(v)}-\frac{d_v'}{\varepsilon_{T_{\pi'}^*}(v)}+
\frac{d_w''}{\varepsilon_{T_{\pi''}}(w)}-\frac{d_w'}{\varepsilon_{T_{\pi'}^*}(w)}+
\frac{d_u''}{\varepsilon_{T_{\pi''}}(u)}-\frac{d_u'}{\varepsilon_{T_{\pi'}^*}(u)}\\
&\geq\frac{d_v'-1}{\varepsilon_{T_{\pi'}^*}(v)}-\frac{d_v'}{\varepsilon_{T_{\pi'}^*}(v)}+
\frac{1}{\varepsilon_{T_{\pi'}^*}(w)}-\frac{1}{\varepsilon_{T_{\pi'}^*}(w)}+
\frac{d_u'+1}{\varepsilon_{T_{\pi'}^*}(u)}-\frac{d_u'}{\varepsilon_{T_{\pi'}^*}(u)}\\
&=\frac{1}{\varepsilon_{T_{\pi'}^*}(u)}-\frac{1}{\varepsilon_{T_{\pi'}^*}(v)}.
\end{array}$$
By the proof of Theorem~\ref{greedy-t}, we can see $\varepsilon_{T_{\pi'}^*}(u)\leq\varepsilon_{T_{\pi'}^*}(v)$. So, $\xi^{ce}(T_{\pi''})\geq \xi^{ce}(T_{\pi'}^*)$.

Hence, $ \xi^{ce}(T_{\pi''}^*) \geq \xi^{ce}(T_{\pi''}) \geq \xi^{ce}(T_{\pi'}^*).$ \ezm

\section{The connective eccentricity index of trees with a given number of branching vertices}

A vertex of a tree $T$ with degree $3$ or greater is called a branching vertex of $T$. For such a tree $T$, it is easy to find that $r\leq \frac{n}{2}-1$. Note that each tree different from the path possesses at least one branching vertices. In the following, we will find a lower bound and an upper bound for
the connective eccentricity index of an $n$-vertex tree with a given number of branching vertices.

Let $\mathcal{BT}_{n,r}$ be the set of all $n$-vertex trees with
exactly $r$ branching vertices. $F(n,r)$ is the greedy caterpillar with degree sequence $\overline{d}=(\overset{r}{\overbrace{3,\cdots,3}},2,\cdots,2,1,\cdots,1)$, and
$B(n,r)$ is the greedy tree with degree  sequence $\overline{d}=(\overset{r}{\overbrace{n-2r+1,3,\cdots,3}}, 1,\cdots,1)$, see Figure 7. Clearly, $F(n,r), B(n,r)\in\mathcal{BT}_{n,r}$ and $B(n,1)=S_n$.

\begin{thm}\label{thm4.1}
If $T\in \mathcal{BT}_{n,r}$ and $1\leq r \leq \frac{n}{2}-1$, then
 $$ \xi^{ce}(T)\geq \xi^{ce}(F(n,r))$$
with equality if and only if $T=F(n,r)$.
\end{thm}
\zm Let $T\in \mathcal{BT}_{n,r}$ be a tree with the maximal connective eccentricity  index. $P=v_0v_1\cdots v_t$ is a longest path
in $T$, and $u_1,u_2,\cdots,u_r$ are all branching vertices of $T$.

First, we show that $d(v)\leq 3$ for $u\in V(T)$. If there is a vertex $u_i$ with $d(u_i)>3$ and $w$ is its neighbor and $w\not\in P$ (See Figure 8). Create a new tree $T'$ (See Figure 8) from $T$ by replacing the edge $u_iw$ with $v_tw$. Notice that $T$ and $T'$ have the same number of branch vertices, and $\varepsilon_{T'}(s) \geq \varepsilon_{T}(s)$ for any vertex $s \in V$ since $P$ is a longest path in $T$. For any vertex $s \in V-{u_i}$, $d_{T'}(s) = d_{T}(s)$ and $d_{T'}(u_i) = d_{T}(u_i)-1$. So, we have
$$\begin{array}{ll}
\xi^{ce}(T) - \xi^{ce}(T') &\geq \frac{d_T(u_i)}{\varepsilon_{T}(u_i)}-\frac{d_{T'}(u_i)}{\varepsilon_{T'}(u_i)}+ \frac{d_T(w)}{\varepsilon_{T}(w)}-\frac{d_{T'}(w)}{\varepsilon_{T'}(w)}\\
&\geq \frac{d_T(u_i)}{\varepsilon_{T'}(u_i)}-\frac{d_{T'}(u_i)}{\varepsilon_{T'}(u_i)}+ \frac{d_T(w)}{\varepsilon_{T'}(w)}-\frac{d_{T'}(w)}{\varepsilon_{T'}(w)}\\
&=\frac{1}{\varepsilon_{T'}(u_i)}>0
\end{array}$$
a contradiction to the extremal property of $T$.

From above, we know that $T$ is a tree with the degree sequence $\overline{d}=(\overset{r}{\overbrace{3,\cdots,3}},2,\cdots,\\ 2,1,\cdots,1)$.
By Theorem \ref{thm3-2}, we have the greedy caterpillar with the degree sequence $\overline{d}=(\overset{r}{\overbrace{3,\cdots,3}},2,\cdots,2,1,\cdots,1)$.
The result is true. \ezm

\begin{thm}\label{thm4.2}
If $T\in\mathcal{BT}_{n,r}$ and $1\leq r \leq \frac{n}{2}-1$, then
$$ \xi^{ce}(T)\leq \xi^{ce}(B(n,r)).$$
\end{thm}
\zm Let $T\in\mathcal{BT}_{n,r}$ be a tree with the minimal connective eccentricity index. Note that every pendant path in $T$ is
a pendant edge by Lemma \ref{lem2.2}.

We first show that $T$ has no vertex of degree two. If $v$ is a vertex of degree two in $T$, then there is a branching
vertex $u$ in $T$ such that $\varepsilon(u)>\varepsilon(v)$ and its
neighbors except one are pendant vertices $v_1,\cdots,v_k$, where
$k=deg(u)-1$ (see Figure 9). Create a new tree $T'$ from $T$ by replacing edges $uv_i (1\leq i\leq k)$ with $vu_i(1\leq i\leq k)$. Notice that $T'\in\mathcal{BT}_{n,r}$ with $\varepsilon_{T'}(s)\leq \varepsilon_{T}(s)$ for any vertex $s \in V(T)$, $d_{T'}(s) = d_{T}(s)$ for any vertex $s \in V-{u,v}$ and $d_{T'}(u) = d_{T}(u)-k=1$, $d_{T'}(v) = d_{T}(v)+k=2+k$, $\varepsilon_{T'}(u)\leq \varepsilon_{T'}(v)$.  So, we have
$$\begin{array}{ll}
\xi^{ce}(T) - \xi^{ce}(T') &\leq \frac{d_T(u)}{\varepsilon_{T}(u)}-\frac{d_{T'}(u)}{\varepsilon_{T'}(u)}+ \frac{d_T(v)}{\varepsilon_{T}(v)}-\frac{d_{T'}(v)}{\varepsilon_{T'}(v)}\\
&\leq \frac{d_T(u)}{\varepsilon_{T'}(u)}-\frac{d_{T'}(u)}{\varepsilon_{T'}(u)}+ \frac{d_T(v)}{\varepsilon_{T'}(v)}-\frac{d_{T'}(v)}{\varepsilon_{T'}(v)}\\
&=k(\frac{1}{\varepsilon_{T'}(u)}-\frac{1}{\varepsilon_{T'}(v)})<0
\end{array}$$
a contradiction to the extremal property of $T$.

From above, we know that $T$ is a tree with degree sequence $\overline{d}=(d_1,\cdots,d_r,1,\cdots,1)$.
By Theorem \ref{thm3-3} and Theorem \ref{greedy-t}, we have the greedy tree with the degree sequence
$\overline{d}=(\overset{r}{\overbrace{n-2r+1,3,\cdots,3}},1,\cdots,1)$.
The result holds. \ezm

\begin{figure}[ht!]
\begin{center}
\includegraphics[width=6cm]{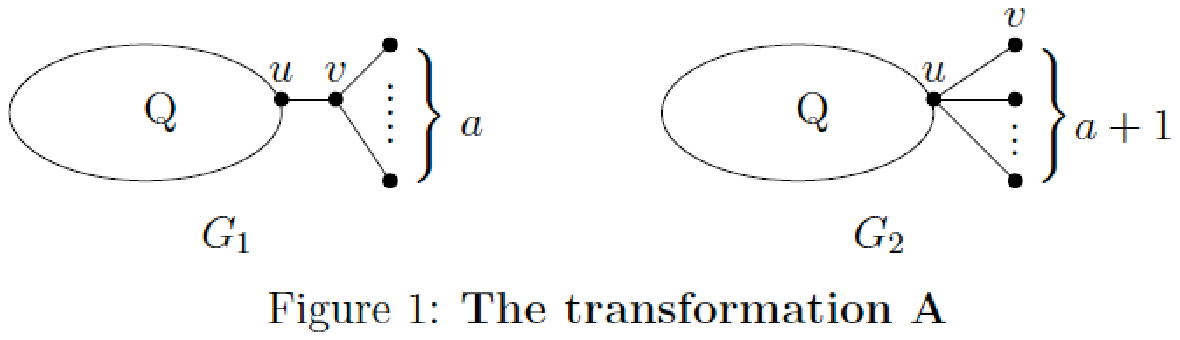}
\includegraphics[width=6cm]{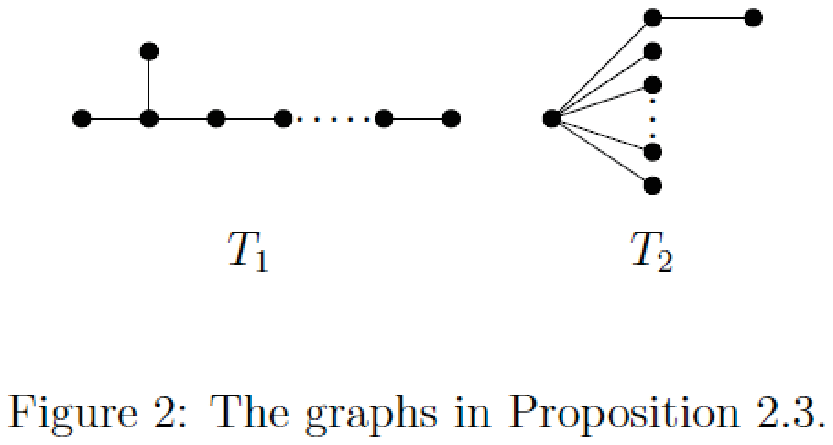}
\includegraphics[width=6cm]{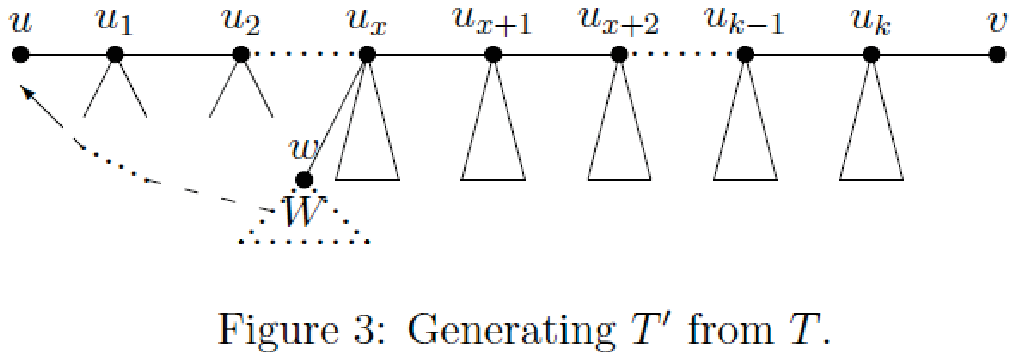}
\includegraphics[width=6cm]{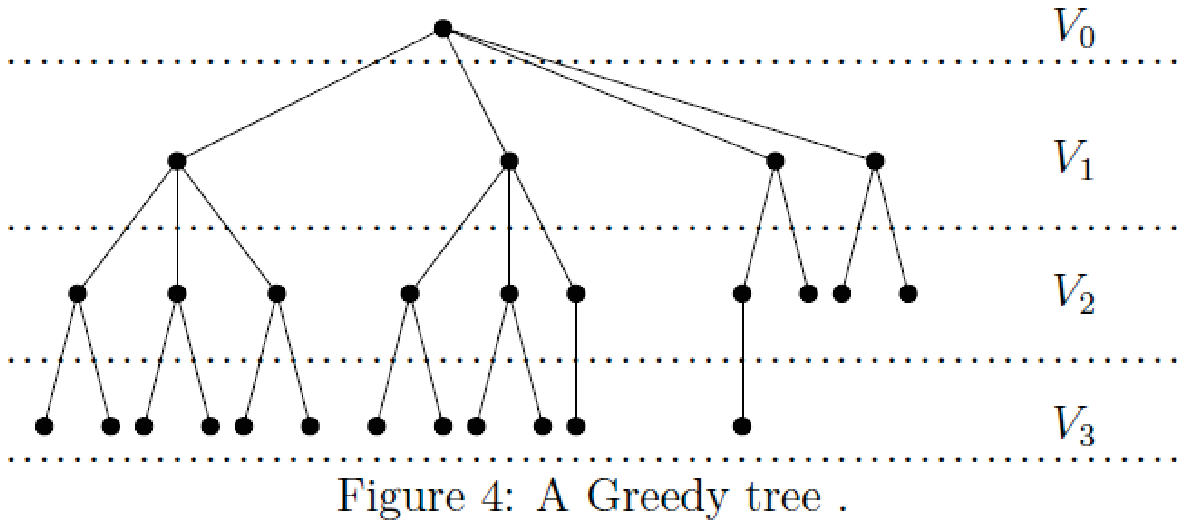}
\includegraphics[width=6cm]{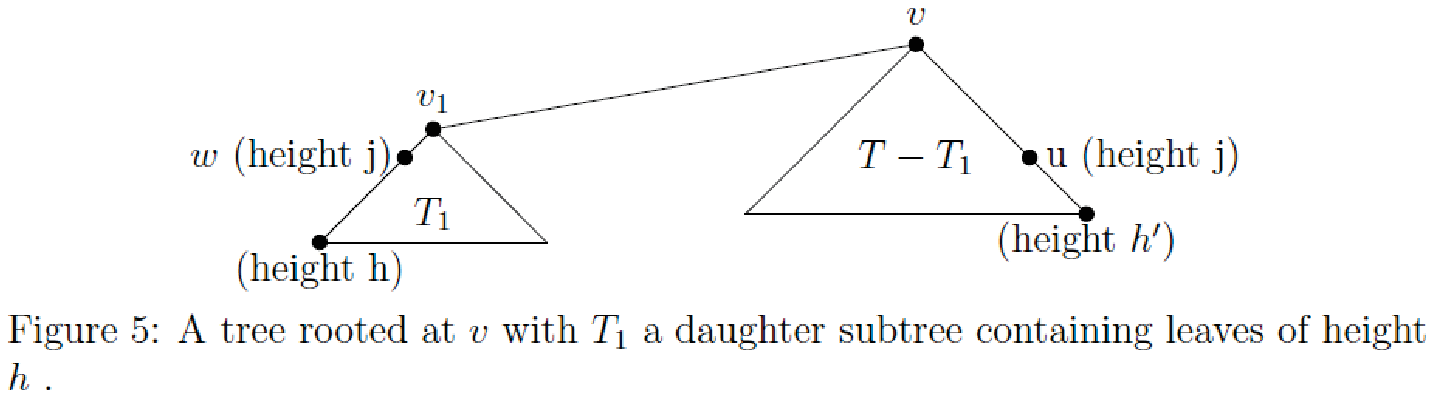}
\includegraphics[width=6cm]{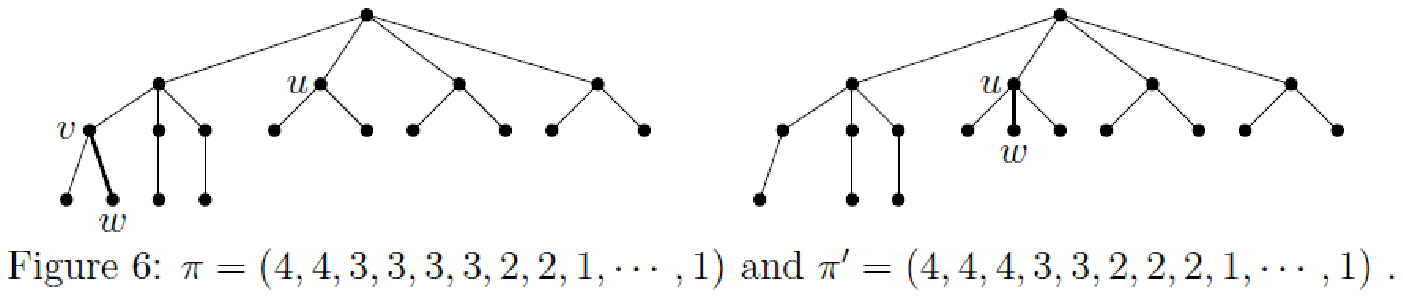}
\includegraphics[width=6cm]{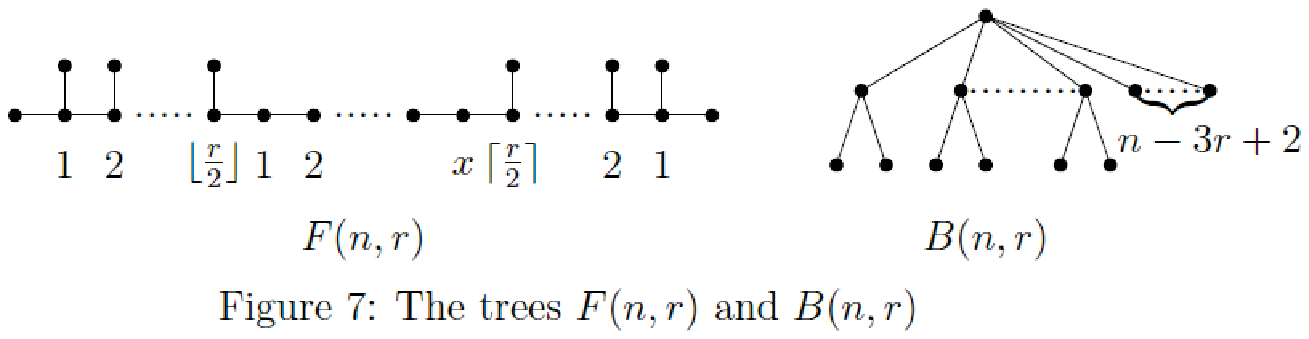}
\includegraphics[width=6cm]{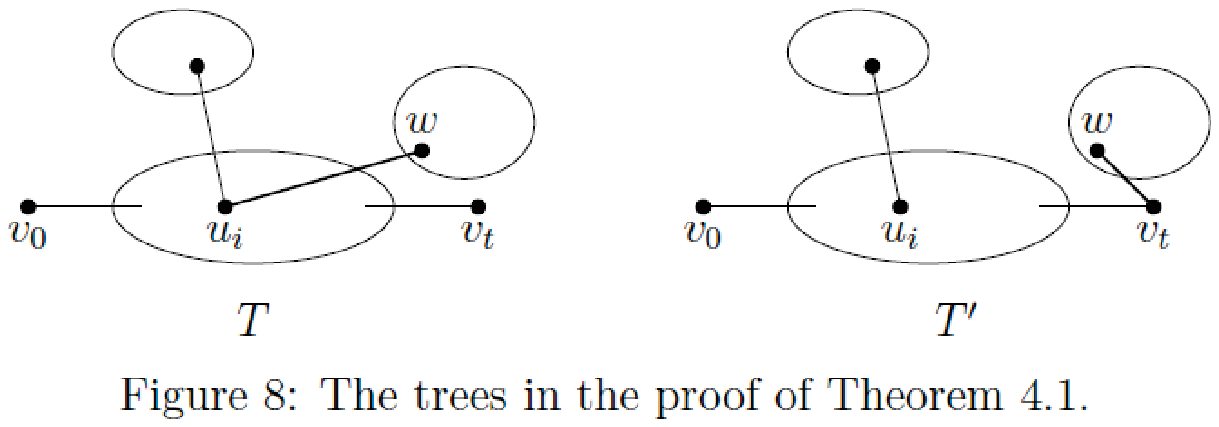}
\includegraphics[width=6cm]{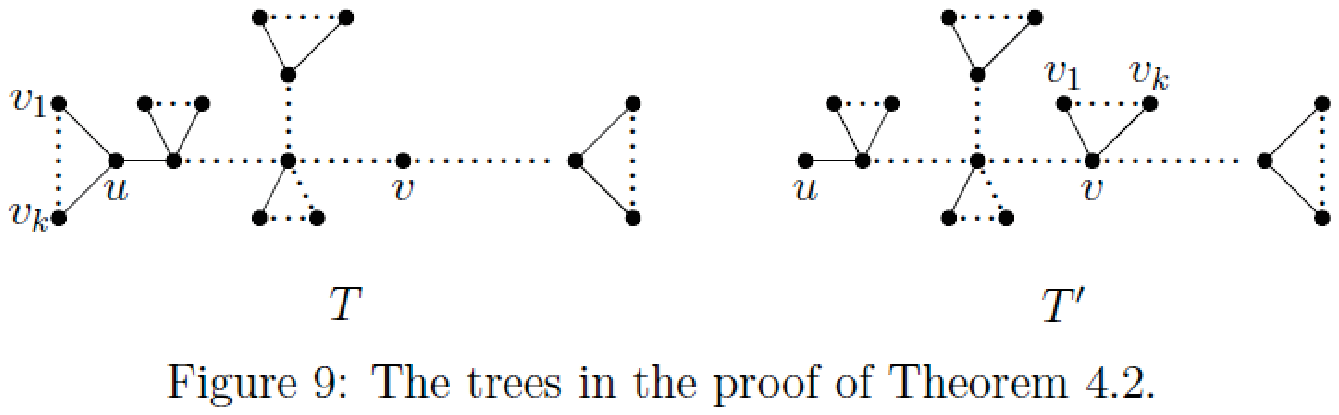}
\end{center}
\end{figure}

\end{document}